\documentclass[letterpaper]{amsart}
\usepackage{amsxtra}
\usepackage{mathtools}
\usepackage{extarrows}
\usepackage{upgreek}
\usepackage{color}
\usepackage{xcolor}
\usepackage{enumitem}
\usepackage{subcaption}
\usepackage{rotating}
\setlist[enumerate]{label=\upshape(\arabic*)}
\allowdisplaybreaks[4]
\usepackage[colorlinks=true, citecolor=blue]{hyperref}

\usepackage{amsthm}
\usepackage{amssymb}
\usepackage{amsmath}
\numberwithin{equation}{section}
\usepackage{upref}
\usepackage{tikz}
\usetikzlibrary{calc}

\newtheorem{theorem}{Theorem}[section]
\newtheorem{theorem-definition}[theorem]{Theorem-Definition}
\newtheorem{lemma}[theorem]{Lemma}
\newtheorem{proposition}[theorem]{Proposition}
\newtheorem{corollary}[theorem]{Corollary}

\newtheorem{question}[theorem]{Question}

\newtheorem{definition}[theorem]{Definition}
\newtheorem{eg}[theorem]{Example}
\newtheorem{remark}[theorem]{Remark}

\begin{document}

\title[On Bott--Samelson rings for Coxeter groups
]{On Bott--Samelson rings for Coxeter groups
}



\author{Tao Gui}
\address{(Tao Gui) \newline \indent Beijing International Center for Mathematical Research, Peking University, No.\ 5 Yiheyuan Road, Haidian District, Beijing 100871, P.R. China}
\email{guitao(at)amss(dot)ac(dot)cn}

\author{Lin Sun}
\address{(Lin Sun) \newline \indent Fudan University, 220 Handan Road, Yangpu District, Shanghai 200433, P.R. China}
\email{lsun21(at)m(dot)fudan(dot)edu(dot)cn}

\author{Shihao Wang}
\address{(Shihao Wang) \newline \indent Tsinghua University, No.\ 30 Haidian District, Beijing 100084, P.R. China}
\email{wangshih21(at)mails(dot)tsinghua(dot)edu(dot)cn}

\author{Haoyu Zhu}
\address{(Haoyu Zhu) \newline \indent School of Mathematics (Zhuhai), Sun Yat-sen University, Zhuhai 519082, Guangdong,
P.R. China}
\email{zhuhy86(at)mails(dot)sysu(dot)edu(dot)cn}

\begin{abstract}
We study the cohomology ring of the Bott--Samelson variety. We compute an explicit presentation of this ring via Soergel's result, which implies that it is a purely combinatorial invariant. We use the presentation to introduce the Bott--Samelson ring associated with a word in an arbitrary Coxeter system by generators and relations. In general, it is a split quadratic complete intersection algebra with a triangular pattern of relations. By a result of Tate, it follows that it is a Koszul algebra and we provide a quadratic (reduced) Gr{\"o}bner basis. Furthermore, we prove that it satisfies the whole K\"ahler package, including the Poincar\'e duality, the hard Lefschetz theorem, and the Hodge--Riemann bilinear relations.
\end{abstract}

\maketitle

\setcounter{tocdepth}{2}

\section{Introduction}
Bott--Samelson varieties are smooth projective varieties introduced firstly by Bott and Samelson \cite{bott1955cohomology,bott1958applications} as smooth manifolds, which are useful in studying the topology of compact Lie groups and symmetric spaces. Hansen \cite{hansen1973cycles} and Demazure \cite{demazure1974desingularisation} considered them as algebraic varieties and used them to construct resolutions of Schubert varieties, which are certain (usually) singular subvarieties of flag varieties and are important in the Schubert calculus and the representation theory of Lie-theoretic objects.

Motivated by the Kazhdan--Lusztig positivity conjecture on the Kazhdan--Lusztig polynomials and by seeking an algebraic proof of the celebrated Kazhdan--Lusztig conjecture on the characters of simple highest weight modules of complex semi-simple Lie algebras \cite{kazhdan1979representations}, Soergel \cite{soergel1990kategorie,soergel1992combinatorics,zbMATH01131858,soergel2007kazhdan} in the 1990s and 2000s developed his theory what is known today as Soergel bimodules. The theory gives a purely algebraic/combinatorial description of the (equivariant) intersection cohomology of Schubert varieties and works for arbitrary Coxeter groups with (currently) no geometric interpretation. 

The starting point of the theory of Soergel bimodules based on an algebraic description of the cohomology of the Bott--Samelson varieties (see \eqref{def-ten} and Remark \ref{rmk-Soergel}) and modeled on the fundamental decomposition theorem of Beilinson, Bernstein, Deligne and Gabber \cite{beilinson2018faisceaux}. Elias and Williamson \cite{elias2014hodge} proved an algebraic version of the decomposition theorem for Soergel bimodules by establishing remarkable Hodge theoretic properties of Soergel bimodules, hence proving the Kazhdan--Lusztig conjectures and completing the program initiated by Soergel.
See \cite{elias2016kazhdan} for a good survey. In addition, they presented the category of Bott--Samelson bimodules by generators and relations using planar diagrammatics, which, by taking the Karoubi envelope, provides a presentation of the category of Soergel bimodules \cite{elias2016soergel}. This approach is very useful for computations. 

In this paper, we study the \emph{ring structure}
of the Bott--Samelson module for arbitrary Coxeter groups. We compute an explicit presentation of this ring via Soergel's result (see Section \ref{sec-pre-Tensor} for the definition of this ring as iterated tensor products by Soergel). 

\begin{theorem} 
For any expression $\underline{w}=$ $\left(s_1, \ldots, s_n\right)$ in an arbitrary Coxeter system $(W, S)$, the Bott--Samelson ring $\overline{\mathrm{BS}}(\underline{w})$, defined in \eqref{def-ten}, has the following quadratic presentation:
\begin{equation} 
\overline{\mathrm{BS}}\left(s_1, s_2, \ldots, s_n\right) \cong \frac{\mathbb{R}\left[x_1, x_2, \ldots, x_n\right]}{\left(q_1, q_2, \ldots, q_n\right)},
\end{equation}
where $\{q_m \mid m=1,2, \ldots, n\}$ is given by
\begin{equation} \label{eq-rel}
\begin{split}
    q_m& := 
    x_m^2-\left(\sum_{l=1}^{m-1}q_{lm}x_l\right)^2=x_m^2-\sum_{i, j=1}^{m-1}(2-\delta_{ij})q_{im}q_{jm}x_i x_j.
\end{split}
\end{equation}
Here $q_{lm}$ is defined in \eqref{eq-quad} and $\delta_{ij}$ is the Kronecker delta.
\end{theorem}

The above theorem implies the Bott--Samelson ring is a purely combinatorial invariant (see Corollary \ref{cor-com}), which is not obvious from the cohomological interpretation nor Soergel's iterated tensor product presentation. See an open question (Question \ref{ques-iso}) related to this corollary. We use the above theorem to introduce the Bott--Samelson ring $\overline{\mathrm{BS}}(\underline{w})$ associated with a word $\underline{w}$ in an arbitrary Coxeter system by generators and relations, see section \ref{sec-presen}. In general, it can be seen from the above presentation that it is a split\footnote{This means that the generating polynomials of the defining ideal split into products of linear factors.} quadratic complete intersection algebra with a triangular pattern of relations, see section \ref{sec-Koszul}. By a result of Tate, it follows that it is a Koszul algebra (Theorem \ref{thm-Koszul}) and we provide a quadratic (reduced) Gr{\"o}bner basis for the defining ideal of relations (Proposition \ref{prop-Grobner}). Furthermore, we show that these rings ``behave as the cohomology of a Bott--Samelson variety'' even though no Bott--Samelson variety exists for most (more precisely, non-crystallographic, see Remark \ref{rmk-cry}) Coxeter systems, by proving that they satisfy the whole K\"ahler package of Hodge theory, including the Poincar\'e duality, the hard Lefschetz theorem, and the Hodge--Riemann bilinear relations.

\begin{theorem} \label{thm-Hodge}
     For any expression $\underline{w}=$ $\left(s_1, \ldots, s_n\right)$ in an arbitrary Coxeter system $(W, S)$, there exists a sequence $\{c_i\}_{i = 1}^n$ with $c_i>0$, such that the following Kähler package holds  with respect to $\ell_n=\sum_{1\leq i\leq n} c_i x_i \in \overline{\mathrm{BS}}^1(\underline{w})$.
     \begin{enumerate}
         \item (Poincaré duality theorem) For every non-negative $k \leq n / 2$, the bilinear pairing defined by
$$\overline{\mathrm{BS}}^k(\underline{w}) \times \overline{\mathrm{BS}}^{n-k}(\underline{w}) \longrightarrow \mathbb{R}, \quad\left(\eta_1, \eta_2\right) \longmapsto \operatorname{deg}\left(\eta_1 \eta_2\right)$$
         is non-degenerate, where the degree isomorphism $\operatorname{deg}: \overline{\mathrm{BS}}^n(\underline{w}) \rightarrow \mathbb{R}$ is defined by taking the coefficient of the monomial $x_1 x_2 \cdots x_n$ \footnote{In geometric setting, this bilinear pairing coincides with the Poincaré pairing, defined by taking cup product of two cohomology classes and evaluating on the fundamental class.}.
         \item (Hard Lefschetz theorem) For every non-negative $k \leq n / 2$, the multiplication map
$$
\overline{\mathrm{BS}}^k(\underline{w}) \longrightarrow \overline{\mathrm{BS}}^{n-k}(\underline{w}), \quad \eta \longmapsto \ell_n^{n-2 k} \eta
$$
is an isomorphism.
\item (Hodge--Riemann relations) For every non-negative $k \leq n / 2$, the bilinear form defined by
$$
\overline{\mathrm{BS}}^k(\underline{w})  \times \overline{\mathrm{BS}}^k(\underline{w})  \longrightarrow \mathbb{R}, \quad\left(\eta_1, \eta_2\right) \longmapsto \operatorname{deg}\left(\ell_n^{n-2 k} \eta_1 \eta_2\right)
$$
is $(-1)^k$-definite when restricted to the primitive subspace 
$$
P_{\ell_n}^{k}=\operatorname{ker}\left(\ell_n^{n-2k+1}\right) \subset\overline{\mathrm{BS}}^k(\underline{w}) .
$$
     \end{enumerate}
\end{theorem}

The Poincar\'e duality follows from the orthogonality of a graded basis of $
\overline{\mathrm{BS}}(\underline{w})$, while we prove the hard Lefschetz and Hodge--Riemann relations together as a single package by induction on the length of the word $\underline{w}$ based on an refined form of \cite[Lemma 5.2]{elias2014hodge} and our presentation.

After completing our work, we became aware of the paper \cite{zbMATH07303202}. In their study, the authors considered the Bott--Samelson ring for a finite reflection group (generated by reflections of finite order) over a field. They computed a presentation of the ring and showed that it has the strong Lefschetz property. Their motivation stemmed from questions regarding embeddings of complete intersections and the inheritance of the strong Lefschetz property. In contrast, the aim of the present paper is to study the Bott--Samelson ring for general Coxeter groups, which, of course, includes infinite ones. By focusing on this setting, we achieve a more uniform and concrete presentation determined by the Coxeter data, which implies the combinatorial invariance. Additionally, we prove the Koszul property, provide a quadratic (reduced) Gr{\"o}bner basis, and establish the whole K\"ahler package, especially the Hodge--Riemann bilinear relations, which were not addressed in their paper.

One of the motivations of this paper is the hope of the first author for understanding and explaining the parallels between the world of Coxeter groups and the world of matroids discovered recently by several people. Given a (complex) hyperplane arrangement, the (now-called) ``matroid Schubert variety'' is a natural compactification of the hyperplane arrangement complement, which is usually singular. De Concini and Procesi \cite{de1995wonderful} constructed canonical resolutions of matroid Schubert varieties called the wonderful models. Feichtner and Yuzvinsky \cite{feichtner2004chow} computed a presentation of the cohomology ring of the wonderful model,  which depends only on the underlying combinatorial (that is, ``matroid'') structure, thus developing a combinatorial cohomology theory for matroids known today as the Chow ring of a matroid. Adiprasito, Huh, and Katz \cite{adiprasito2018hodge} proved that the Chow ring of any matroid satisfies the Hodge theory, regardless of whether it has a geometric origin (that is, comes from a hyperplane arrangement) or not, thus developing a ``smooth'' combinatorial Hodge theory for matroids. They used the Hodge--Riemann relations to resolve the Heron--Rota--Welsh conjecture for arbitrary matroids. Dotsenko \cite{dotsenko2022homotopy} conjectured that the Chow ring of any matroid is Koszul, which is proven by Matthew Mastroeni and Jason McCullough \cite{mastroeni2023chow}. This paper can be seen as to exploit a parallel  combinatorial cohomology theory and a ``smooth'' combinatorial Hodge theory for arbitrary Coxeter groups. We hope that this paper could shed some light on the comparison of the two worlds and might be helpful in pursuing some benefits from each other, see Question \ref{ques-brick!} for more directions on this perspective.

The remainder of the paper is organized as follows. In Section \ref{sec-preli}, we introduce the notations and some preliminary results. In Section \ref{sec-presen}, we compute the presentations of Bott--Samelson rings. In Section \ref{sec-Koszul}, we prove that Bott--Samelson rings are Koszul and provide a quadratic Gr{\"o}bner basis for the defining ideal of relations. In Section \ref{sec-Hodge}, we prove that Bott--Samelson rings satisfy K\"ahler package. In Section \ref{sec-question}, we present some future directions and state some open questions.

\subsection*{Acknowledgments}
The first author would like to express his gratitude to Geordie Williamson for valuable communications and would like to thank Ben Elias for helpful discussions. This research was carried out as part of the PACE program in the summer of 2024 at Peking
University, Beijing. We thank Yibo Gao who organized the program and the university for providing a nurturing
academic atmosphere that significantly contributed to the progression of the research.

\section{Preliminaries} \label{sec-preli}

In this section, we collect our notations and some preliminary results for later use. A good  reference for this section is \cite[Part I]{elias2020introduction}.

\subsection{Coxeter Systems and the Geometric Representation} \label{sec-pre-Coxeter}
A \emph{Coxeter system} $(W, S)$ is a group $W$ with a finite set $S \subset W$ of generators of $W$, such that
$$
\left.W=\langle s \in S \quad | \quad (s t)^{m_{s t}}=\mathrm{id} \text { for any } s, t \in S \text { with } m_{s t}<\infty\right\rangle,
$$
where $m_{s s}=1$ for each $s \in S$, and $m_{s t}=m_{t s} \in\{2,3, \ldots\} \cup\{\infty\}$ for $s \neq t \in S$. It follows that $m_{s t}$ is precisely the order of the element $s t$, and when $m_{s t}=\infty$, there is no corresponding relation between $s$ and $t$. The elements of $S$ are called \emph{simple reflections}. For each $w \in W$, one can write $w=s_1 \cdots s_k$ for some $s_1, \ldots, s_k \in S$. The sequence $\left(s_1, \ldots, s_k\right)$ is called an \emph{expression} or a \emph{word} for $w$. We use a notational shorthand $\underline{w}$ to denote the sequence $\left(s_1, \ldots, s_k\right)$, when the product $s_1 \cdots s_k$ is equal to $w$. The \emph{length} of $w$, denoted by $\ell(w)$, is the minimal $k$ for which $w$ admits an expression $\left(s_1, \ldots, s_k\right)$. Any expression for $w$ with this minimal length $\ell(w)$ is called a \emph{reduced expression}. In particular, $\ell(w)=0$ if and only if $w=\mathrm{id}$.

Let $V$ be the real vector space with a basis $\left\{\alpha_s \mid s \in S\right\}$ indexed by the finite set $S$ of simple reflections of a Coxeter system $(W, S)$. These basis elements $\alpha_s$ are known as \emph{simple roots}. Equip $V$ with the symmetric bilinear form $(-,-)$ determined by
$$
\left(\alpha_s, \alpha_t\right)=-\cos \frac{\pi}{m_{s t}}.
$$
When $m_{s t}=\infty$, we use the convention that $\frac{\pi}{m_{s t}}=0$, so that in this case $\left(\alpha_s, \alpha_t\right)=-1$. Note that $\left(\alpha_s, \alpha_s\right)=1$. The \emph{geometric representation} of the Coxeter system $(W, S)$ is the representation $V$ of $W$ defined to be the action of $W$ on $V$, where each simple reflection $s \in S$ acts by reflection along $\alpha_s$. That is,
$$
s(\lambda)=\lambda-2\left(\lambda, \alpha_s\right) \alpha_s , \quad \forall \lambda \in V.
$$
It follows that the symmetric bilinear form $(-,-)$ on $V$ is $W$-invariant for the geometric representation. It is known that for any Coxeter system $(W, S)$, the geometric representation is faithful, and $W$ is finite if and only if the bilinear form $(-,-)$  is positive definite.

\subsection{Invariant polynomials and Demazure operators} \label{sec-pre-Bott}

For a $\mathbb{Z}$-graded vector space $M:=\bigoplus_{i \in \mathbb{Z}} M^i$, we use the notation $M(d)$ for the graded vector space with graded pieces $M(i)^j:=M^{i+j}$. \footnote{That is, if one imagine the graded pieces of $M$ arranged vertically, then $M(1)$ is obtained by shifting $M$ down.}

Fix a Coxeter system $(W, S)$, let 
$$
R=\operatorname{Sym}(V)=\bigoplus_{i \in \mathbb{Z}_{\geq 0}} \operatorname{Sym}^i(V)
$$
be the symmetric algebra of the geometric representation $V$. We view it as a graded algebra with the convention\footnote{Note that our convention is different with the one used in \cite{elias2020introduction}.} that $\operatorname{deg}(V)=1$. In other words, $R$ is the polynomial ring
$$
R=\mathbb{R}\left[\alpha_s \mid s \in S\right]
$$
with grading
$$
\operatorname{deg} \alpha_s=1
$$
for all $s \in S$. The $W$-action on $V$ induces a $W$-action on $R$ via
$$
w \cdot \prod_{s \in S} \alpha_s^{k_s}=\prod_{s \in S}\left(w\left(\alpha_s\right)\right)^{k_s}
$$
on monomials, and then extended linearly to polynomials.

Let $R^W$ be the subring of $W$-invariant polynomials in $R$ and let $R^s$ be the subring of $s$-invariant polynomials in $R$. It is not difficult to show that for every $s \in S, R^s$ is generated by $\alpha_s^2$ and the elements 
$$
\alpha_t+\left(\cos \frac{\pi}{m_{s t}}\right) \alpha_s \quad \text { for all } \quad t \in S \backslash\{s\}.
$$

For $s \in S$, the \emph{Demazure operator} $\partial_s$ (also called the \emph{BGG operator}, or the \emph{divided difference operator}) is the graded map
$$
\begin{aligned}
\partial_s: R & \rightarrow R^s(-1), \\
f & \mapsto \frac{f-s(f)}{\alpha_s} .
\end{aligned}
$$
It is not difficult to see that $\partial_s$ is well-defined.

The Demazure operator is useful in constructing projection operators in a decomposition of $R$ as an graded $R^s$-module into a direct sum $R \simeq R^s \oplus$ $R^s \cdot \alpha_s$ into $s$-invariants and $s$-anti-invariants, given by 
\begin{equation} \label{eq-decom}
  \begin{aligned}
f&=\frac{f+s(f)}{2}+\frac{f-s(f)}{2}\\
&=\partial_s\left(f \frac{\alpha_s}{2}\right)+\frac{\alpha_s}{2} \partial_s(f).
\end{aligned}  
\end{equation}

We collect here some important properties of Demazure operators that will be used later.

\begin{lemma} \label{lem-Demazure}
Let $s \in S$. Then $\partial_s$ satisfies the following properties.
\begin{enumerate}
    \item $\partial_s\left(\alpha_t\right)=-2\cos \frac{\pi}{m_{s t}}$ for any $t \in S$. In particular, $\partial_s\left(\alpha_s\right)=2$.
    \item Let $\operatorname{p}_s \colon R \simeq R^s \oplus R^s\cdot \alpha_s \to R^s$ be the projection map, that is, for $f \in R$ we have $\operatorname{p}_s(f) = \frac{f+s(f)}{2}$. Then $\operatorname{p}_s(\alpha_t) = \alpha_t + \cos \frac{\pi}{m_{st}} \alpha_s$.
    \item  $\partial_s$ is an $R^s$-bimodule map.
\end{enumerate}
\end{lemma}

\begin{proof}
    One directly checks that\[\partial_s(\alpha_t) = \frac{\alpha_t-(\alpha_t+2\cos \frac{\pi}{m_{st}}\alpha_s)}{\alpha_s} = -2\cos \frac{\pi}{m_{st}}\]
    and
    \[\operatorname{p}_s(\alpha_t) = \frac{\alpha_t+(\alpha_t+2\cos \frac{\pi}{m_{st}}\alpha_s)}{2} = \alpha_t+\cos \frac{\pi}{m_{st}}\alpha_s.\]

    For $f \in R$ and $g \in R^s$, we have $\partial_s(gf) = \frac{gf-s(gf)}{2} = \frac{gf-g\cdot s(f)}{2} = g\partial_s(f)$ and similarly $\partial_s(fg) = \partial_s(f)g$, which show that $\partial_s$ is an $R^s$-bimodule map.
\end{proof}

\subsection{Bott--Samelson rings as iterated tensor products} \label{sec-pre-Tensor}

Given an expression $\underline{w}=$ $\left(s_1, \ldots, s_n\right)$ in a Coxeter system $(W, S)$, the corresponding \emph{Bott--Samelson ring}, denoted by $\overline{\mathrm{BS}}(\underline{w})$, is the graded $\mathbb{R}$-algebra given by
\begin{equation} \label{def-ten}
\begin{aligned}
  \overline{\mathrm{BS}}(\underline{w}) :=& \mathbb{R} \otimes_{R} R 
 \otimes_{R^{s_1}} R \otimes_{R^{s_2}} \cdots \otimes_{R^{s_n}} R\\
  =& \mathbb{R} \otimes_{R^{s_1}} R \otimes_{R^{s_2}} \cdots \otimes_{R^{s_n}} R,  
\end{aligned}
\end{equation}
where the $R$-module structure on $\mathbb{R}$ is given by the natural quotient map 
$R \rightarrow \mathbb{R}$.

The multiplicative structure is defined component-wise by 
$$
\left(a \otimes f_1 \otimes \cdots \otimes f_n\right) \cdot\left(b  \otimes g_1 \otimes \cdots \otimes g_n\right)=\left(ab \right) \otimes \left( f_1 g_1\right) \otimes \cdots \otimes\left(f_n g_n \right),
$$
where $a, b \in \mathbb{R}$ and $f_1, \dots, f_n, g_1, \dots, g_n \in R$.

\begin{remark} \label{rmk-cry}
    When $(W, S)$ is a \emph{crystallographic} Coxeter system, that is, $m_{s t} \in$ $\{2,3,4,6, \infty\}$ for all $s \neq t \in S$ (for example, when $W$ is a Weyl group or affine Weyl group), one can define a representation of $W$ over $\mathbb{Z}$ rather than over $\mathbb{R}$ by modifying the simple roots in the geometric representation with possibly different lengths (for example, simple roots of non-simply-laced type Lie algebras). This integral representation can be realized on the Cartan subalgebra of the corresponding complex semisimple Lie algebra. However, as $W$-representations, this representation is isomorphic to the geometric representation over $\mathbb{R}$. 
    It is not difficult to see that if one uses this representation other than the geometric representaion to define the Bott--Samelson ring, one would get an isomorphic ring. Since the geometric representaion behaves in a uniform way for all Coxeter systems, we stick to the geometric representaion to define the Bott--Samelson ring.
\end{remark}

\begin{remark} \label{rmk-Soergel}
    When $(W, S)$ is a \emph{crystallographic} Coxeter system, a theorem of Soergel says that the Bott--Samelson ring is isomorphic to the cohomology ring of the corresponding Bott--Samelson variety. This observation is the starting point of the beautiful theory of Soergel bimodules \cite{soergel1990kategorie,soergel1992combinatorics,zbMATH01131858,soergel2007kazhdan}.
\end{remark}

\section{A Quadratic Presentation and Generators-Relations Definition of Bott--Samelson Rings for Arbitrary Coxeter Groups} \label{sec-presen}
Given an expression $\underline{w}=$ $\left(s_1, \ldots, s_n\right)$ in a Coxeter system $(W, S)$, define the simple tensor $x_i \in \overline{\mathrm{BS}}(\underline{w})$ as
$$
x_i=1 \otimes 1 \otimes \cdots \otimes \alpha_{s_i} \otimes \cdots \otimes 1,
$$
where $\alpha_{s_i}$ is in the $(i+1)$-th plot. We will give $\overline{\mathrm{BS}}(\underline{w})$ a quadratic presentation in the last of this section, and we will use it as the generators-relations definition of the Bott--Samelson ring $\overline{\mathrm{BS}}(\underline{w})$.

For convenience and intuition, we may use $a \underset{s_1}{|} f_1 \underset{s_2}{|} \cdots \underset{s_n}{|} f_n$ as in \cite{elias2020introduction} instead of $a \otimes f_1 \otimes \cdots \otimes f_n$, with  ``$\underset{s}{|}$'' meaning a ``wall'' allowing only $R^s$-invariant elements to pass through. Also in our case , for $i \in [n]$, we denote $\underset{i}{|}$ to mean $\underset{s_i}{|}$, and similarly $\partial_i := \partial_{s_i}$, $\operatorname{p}_i := \operatorname{p}_{s_i}$ and $r_{ij} := \cos \frac{\pi}{m_{s_is_j}}$.

For a homogeneous polynomial $f \in R$ and $1 \le m \le n$, by repeatedly using the decomposition \eqref{eq-decom} to move $f$ all the way to the left, we get
\[1\underset{1}{|} 1 \underset{2}{|}\cdots 1 \underset{m-1}{|} f \underset{m}{|} 1 \cdots 1 \underset{n}{|} 1 = 1 \cdot f_0\underset{1}{|} 1 \underset{2}{|} 1 \cdots +\sum_{k = 1}^{m-1} \left(\sum_{1 \leq i_1<\ldots < i_k \leq m-1} d_{i_1\cdots i_k}^{m-1}(f) x_{i_1}\cdots x_{i_k}\right),\]
where $f_0 = \operatorname{p}_1 \circ \operatorname{p}_2 \circ \cdots \circ \operatorname{p}_{m-1} (f)$ and \[d_{i_1\cdots i_k}^{m-1}(f) = \operatorname{p}_1 \circ \cdots \circ \operatorname{p}_{i_1-1} \circ \frac{\partial_{i_1}}{2} \circ \operatorname{p}_{i_1+1} \circ \cdots \circ \operatorname{p}_{i_2-1} \circ \frac{\partial_{i_2}}{2} \circ \cdots (f) \in \mathbb{R}.\]
The above equations hold because repeatedly using the decomposition \eqref{eq-decom} to move $f$ all the way to the left creates $2^{m-1}$ terms. In each time of the decompositions, the invariant part $\operatorname{p}_{i}(-)$ cross the wall, while $\frac{\partial_{i}}{2}(-)$ in the anti-invariant part cross the wall leaving $\alpha_i$ behind. The indices $1 \leq i_1<\ldots < i_k \leq m-1$ record the elements of the subset of $[2]^{m-1}$ choosing anti-invariant part.

Note that $1 \cdot f_0 = 0$. If $f$ is of degree $1$, then $d_{i_1\cdots i_k}^{m-1}(f) = 0$ if $k \geq 2$. In this case, we write\[1\underset{1}{|} 1 \underset{2}{|}\cdots 1 \underset{m-1}{|} f \underset{m}{|} 1 \cdots = \sum_{l = 1}^{m-1} d_l^{m-1}(f) x_l.\]
By definition, one directly checks that
$$
d_l^{m-1}(f) = \begin{cases}
    \frac{\partial_{m-1} f}{2}, & l = m-1; \\ \frac{\partial_l}{2}\left(\operatorname{p}_{l+1} \circ \cdots \circ \operatorname{p}_{m-1}(f)\right), & l<m-1.
\end{cases}
$$

We are especially concerned with the case that $f = \alpha_{s_m}$. Actually, we have the following computation.
\begin{lemma}
    Suppose $l<m-1$, then \begin{equation}\label{eq-computation}
    \begin{split}
        & \operatorname{p}_{l+1} \circ \cdots \circ \operatorname{p}_{m-1}(\alpha_{s_m}) \\ & = \sum_{k = l+1}^{m-1} \left(\sum_{s = 0}^{m-1-k} \sum_{k=j_0<j_1 <\ldots<j_s<j_{s+1} = m}\prod_{\alpha = 0}^s r_{j_\alpha j_{\alpha+1}}\right)\alpha_{s_k}+\alpha_{s_m}
    \end{split}
    \end{equation}
    and
    \begin{equation}\label{eq-dlm}
        d_l^{m-1}(\alpha_{s_m}) = -\sum_{s = 0}^{m-1-l}\left( \sum_{l = j_0<j_1<\ldots<j_s<j_{s+1} = m} \prod_{\alpha = 0}^s r_{j_\alpha j_{\alpha +1}}\right).
    \end{equation}
    In addition, \eqref{eq-dlm} holds when $l = m-1$.
\end{lemma}
\begin{proof}
We prove \eqref{eq-computation} by reverse induction on $l$. Firstly, \eqref{eq-computation} holds when $l = m-2$: $\operatorname{p}_{m-1}(\alpha_{s_m}) = \alpha_{s_m}+r_{m-1,m}\alpha_{s_{m-1}}$. Suppose \eqref{eq-computation} holds for $l < m-1$, then using \eqref{lem-Demazure} we have
    \begin{equation*}
        \begin{split}
            & \operatorname{p}_l \circ \cdots \circ \operatorname{p}_{m-1}(\alpha_{s_m})\\ =& \sum_{k = l+1}^{m-1} \left((\sum_{s = 0}^{m-1-k} \sum_{k=j_0<j_1 <\ldots<j_s<j_{s+1} = m}\prod_{\alpha = 0}^s r_{j_\alpha j_{\alpha+1}})(\alpha_{s_k}+r_{lk}\alpha_{s_l})\right) +(\alpha_{s_m}+r_{lm}\alpha_{s_l}) \\ = &\sum_{k = l}^{m-1} \left(\sum_{s = 0}^{m-1-k} \sum_{k=j_0<j_1 <\ldots<j_s<j_{s+1} = m}\prod_{\alpha = 0}^s r_{j_\alpha j_{\alpha+1}}\right)\alpha_{s_k}+\alpha_{s_m},
        \end{split}
    \end{equation*}
    implying \eqref{eq-computation} holds for $l-1$. Hence \eqref{eq-computation} holds by induction.
    
    Clearly \eqref{eq-dlm} holds when $l = m-1$ by Lemma \ref{lem-Demazure}(1). If $l<m-1$, using \eqref{eq-computation} and Lemma \ref{lem-Demazure}, we have
    \begin{equation*}
        \begin{split}
            & d_l^{m-1}(\alpha_{s_m}) = \frac{\partial_l}{2} \circ \operatorname{p}_{l+1} \circ \cdots \circ \operatorname{p}_{m-1}(\alpha_{s_m}) \\ & =  -\sum_{k = l+1}^{m-1} \left(\sum_{s = 0}^{m-1-k} \sum_{k=j_0<j_1 <\ldots<j_s<j_{s+1} = m}r_{lk}\prod_{\alpha = 0}^s r_{j_\alpha j_{\alpha+1}}\right)-r_{lm} \\ & = -\sum_{s = 0}^{m-1-l}\left( \sum_{l = j_0<j_1<\ldots<j_s<j_{s+1} = m} \prod_{\alpha = 0}^s r_{j_\alpha j_{\alpha +1}}\right).
        \end{split}
    \end{equation*}
\end{proof}

An observation goes as follows: consider the matrix $\Tilde{R}_{\underline{w}} = (\Tilde{r}_{ij})_{1 \leq i,j \leq n}$, where $$\Tilde{r}_{ij} := \begin{cases}
    0, &i \geq j ; \\ r_{ij} = \cos \frac{\pi}{m_{s_is_j}}, &i < j.
\end{cases}$$
Then for $s+1 \in [m-l]$, $\sum_{l = j_0<j_1<\ldots<j_s<j_{s+1} = m} \prod_{\alpha = 0}^s r_{j_\alpha j_{\alpha +1}}$ is the $(l, m)$-th entry of $\Tilde{R}_{\underline{w}}^{s+1}$. Therefore $-d_l^{m-1}(\alpha_{s_m})$ is the $(l, m)$-th entry of $\left(\Tilde{R}_{\underline{w}}+\Tilde{R}_{\underline{w}}^2+\ldots+\Tilde{R}_{\underline{w}}^{m-l}\right)$.  
Since $\Tilde{R}_{\underline{w}}$ is an upper-triangular nilpotent matrix with nilpotent index $n$, $-d_l^{m-1}(\alpha_{s_m})$ is the $(l, m)$-th entry of 
\begin{equation} \label{eq-quad}
 Q_{\underline{w}} = (q_{ij})_{1 \leq i,j \leq n} := \Tilde{R}_{\underline{w}}+\Tilde{R}_{\underline{w}}^2+\ldots+\Tilde{R}_{\underline{w}}^n.   
\end{equation}

We note that $\displaystyle |q_{lm}| \leq \sum_{s = 0}^{m-1-l} \binom{m-1-l}{s} = 2^{m-1-l}$. Besides, we have the following proposition.
\begin{proposition} Let $\underline{w}=$ $\left(s_1, \ldots, s_n\right)$ and $1 \le m \le n$, then $x_m^2 = \left(\sum\limits_{l = 1}^{m-1} q_{lm}x_l\right)^2$ in $\overline{\mathrm{BS}}(\underline{w})$. 
\end{proposition}
\begin{proof}
\begin{equation*}
    \begin{split}
    & x_m^2 = 1 \underset{1}{|} 1 \cdots 1 \underset{m-1}{|}\alpha_{s_m}^2 \underset{m}{|} \cdots = \left(1 \underset{1}{|} 1 \cdots 1 \underset{m-1}{|}\alpha_{s_m} \underset{m}{|} \cdots \right)^2 \\ & = \left(\sum_{l = 1}^{m-1}d_l^{m-1}(\alpha_{s_m})x_l\right)^2 = \left(\sum\limits_{l = 1}^{m-1} q_{lm}x_l\right)^2.
    \end{split}
\end{equation*}
\end{proof}

\begin{theorem-definition} \label{thm-def}
For any expression $\underline{w}=$ $\left(s_1, \ldots, s_n\right)$ in an arbitrary Coxeter system $(W, S)$, the Bott--Samelson ring $\overline{\mathrm{BS}}(\underline{w})$, defined in \eqref{def-ten}, has the following quadratic presentation, which is served as the generators-relations definition of the Bott--Samelson ring $\overline{\mathrm{BS}}(\underline{w})$:
\begin{equation} \label{def-gen}
\overline{\mathrm{BS}}\left(s_1, s_2, \ldots, s_n\right) \cong \frac{\mathbb{R}\left[x_1, x_2, \ldots, x_n\right]}{\left(q_1, q_2, \ldots, q_n\right)},
\end{equation}
where $\{q_m \mid m=1,2, \ldots, n\}$ is given by
\begin{equation*} 
\begin{split}
    q_m& :=
    x_m^2-\left(\sum_{l=1}^{m-1}q_{lm}x_l\right)^2=x_m^2-\sum_{i, j=1}^{m-1}(2-\delta_{ij})q_{im}q_{jm}x_i x_j.
\end{split}
\end{equation*}
Here $q_{lm}$ is defined in \eqref{eq-quad} and $\delta_{ij}$ is the Kronecker delta.
\end{theorem-definition}

\begin{proof}
    From the discussions above, we have a well-defined graded algebra homomorphism from $\displaystyle \frac{\mathbb{R}\left[x_1, x_2, \ldots, x_n\right]}{\left(q_1, q_2, \ldots, q_n\right)} \longrightarrow \overline{\mathrm{BS}}\left(s_1, s_2, \ldots, s_n\right)$ sending $x_i$ to the simple tensor $1 \otimes 1 \otimes \cdots \otimes \alpha_{s_i} \otimes \cdots \otimes 1$. Since $R \simeq R^{s_n} \oplus R^{s_n} \cdot \alpha_s$ as graded left $R^{s_n}$-modules, we have the following isomorphism as graded real vector spaces:
    \begin{equation*}
        \begin{split}
            \overline{\mathrm{BS}}\left(s_1, s_2, \ldots, s_n\right) &= \mathbb{R} \otimes_{R} R 
 \otimes_{R^{s_1}} R \otimes_{R^{s_2}} \cdots \otimes_{R^{s_n}} R\\
 &\cong \mathbb{R}\otimes_{R^{s_1}} R \otimes \ldots \otimes_{R^{s_n}} \left(R^{s_n} \oplus R^{s_n} \cdot \alpha_s\right) \\ & \cong \overline{\mathrm{BS}}\left(s_1, s_2, \ldots, s_{n-1}\right) \oplus \overline{\mathrm{BS}}\left(s_1, s_2, \ldots, s_{n-1}\right)(-1),
        \end{split}
    \end{equation*}
    which has dimension $2^n$  and has graded basis $\{1\otimes\alpha_{s_1}^{\varepsilon_1}\otimes\alpha_{s_2}^{\varepsilon_2}\otimes\cdots\otimes \alpha_{s_n}^{\varepsilon_n} \mid \varepsilon_i \in \{0,1\} ,i \in [n]\}$ by induction. Also it is not difficult to see that $\displaystyle \frac{\mathbb{R}\left[x_1, x_2, \ldots, x_n\right]}{\left(q_1, q_2, \ldots, q_n\right)}$ has real dimension $2^n$ as well, with $\{x_1^{\varepsilon_1}x_2^{\varepsilon_2}\cdots x_n^{\varepsilon_n} \mid \varepsilon_i \in \{0,1\} ,i \in [n]\}$ being a graded basis, see also Proposition \ref{prop-Grobner}. Since we have an $\mathbb{R}$-algebra homomorphism between $\mathbb{R}$-algebras of the same real dimension mapping basis to basis bijectively, it is a graded $\mathbb{R}$-algebra isomorphism. 
\end{proof}

\begin{eg} \label{eg-123}
    When $n=1$, $\overline{\mathrm{BS}}\left(s_1\right) \cong \frac{\mathbb{R}\left[x\right]}{\left(x^2\right)}$. When $n=2$, $\overline{\mathrm{BS}}\left(s_1, s_2\right) \cong \frac{\mathbb{R}\left[x_1, x_2\right]}{\left(x_1^2, x_2^2\right)}$.
    When $n=3$, for the expression $\underline{w}=$ $\left(s_1, s_2, s_3\right)$ in a Coxeter system $(W, S)$, the Bott--Samelson ring $\overline{\mathrm{BS}}(\underline{w})$ has the following presentation:
\begin{equation*} 
\overline{\mathrm{BS}}\left(s_1, s_2, s_3\right) \cong \frac{\mathbb{R}\left[x_1, x_2, x_3\right]}{\left(x_1^2, x_2^2, x_3^2-2 \left( \cos \frac{\pi}{m_{s_1s_3}} \cos \frac{\pi}{m_{s_2s_3}} + \cos\frac{\pi}{m_{s_1s_2}} \cos^2 \frac{\pi}{m_{s_2s_3}} \right) x_1x_2\right)}.
\end{equation*}
\end{eg}

The above theorem has the following immediate corollary, which says that the Bott--Samelson ring $\overline{\mathrm{BS}}(\underline{w})$ only depends on $m_{s t}$'s of those simple reflections appearing in the expression $\underline{w}=$ $\left(s_1, \ldots, s_n\right)$. This is not obvious from the iterated tensor products definition of $\overline{\mathrm{BS}}(\underline{w})$, since the geometric representations of two different Coxeter groups have nothing related at all!

\begin{corollary}[Combinatorial invariance]
\label{cor-com}
Let $\underline{w}=$ $\left(s_1, \ldots, s_n\right)$ be an expression in a Coxeter system $(W, S)$, and $\underline{w^\prime}=$ $\left(s_1^\prime, \ldots, s_n^\prime\right)$ be another expression in a possibly different Coxeter system $(W^\prime, S^\prime)$, if the two matrices $M_w:=\left(m_{s_i s_j}\right)_{1 \leq i \leq j \leq n}$ and $M_w^\prime:=\left(m_{s_i^\prime s_j^\prime}\right)_{1 \leq i \leq j \leq n}$ are equal, then we have 
$$\overline{\mathrm{BS}}(\underline{w}) \cong \overline{\mathrm{BS}}(\underline{w^\prime}).$$
\end{corollary}

\begin{remark}
    It is natural to ask if two arbitrary expressions of the same length are related by braid relations (that is, one can apply a sequence of braid relations to obtain one from the other), whether or not the associated Bott--Samelson rings are isomorphic. It is straightforward to verify that if the length of the expressions are at most 3, or if the two expressions involve only two distinct simple reflections $s$ and $t$ appearing alternatively for $m_{st} < \infty$ times related by a single braid relation, then the corresponding Bott--Samelson rings are isomorphic. However, it is not always the case in general. For example, consider the Coxeter system $(S_3,\{s_1,s_2\})$ and two expressions $\underline{w} = (s_1,s_2,s_1,s_2)$ and $\underline{w}' = (s_1,s_1,s_2,s_1)$. Then one checks that \[\overline{\operatorname{BS}}(\underline{w}) \cong \frac{\mathbb{R}[x_1,x_2,x_3,x_4]}{(x_1^2,x_2^2,x_3^2+\frac{3}{4}x_1x_2,x_4^2-\frac{3}{8}x_1x_2+\frac{3}{8}x_1x_3+\frac{3}{4}x_2x_3)}\]and\[\overline{\operatorname{BS}}(\underline{w}') \cong \frac{\mathbb{R}[x_1,x_2,x_3,x_4]}{(x_1^2,x_2^2,x_3^2,x_4^2+\frac{3}{4}x_2x_3)}.\]
    One can show that they are not isomorphic as graded $\mathbb{R}$-algebras and we leave it as an exercise to the interested readers.
\end{remark}

\section{Bott--Samelson Rings are Koszul and a quadratic Gr{\"o}bner basis} \label{sec-Koszul}

In this section we prove that Bott--Samelson rings are Koszul and use the presentation \eqref{def-gen} to construct a quadratic (reduced) Gr{\"o}bner basis for the defining ideal of relations.

We expand the quadratic term in \eqref{eq-rel}:

\begin{equation} \label{eq-expan}
    q_m = x_m^2-\sum_{1\leq i\leq j\leq m-1}c_{i,j}^{m}x_ix_j,
\end{equation}
where $c_{i,j}^{m}=(2-\delta_{ij})q_{im}q_{jm}$. Since $|q_{lm}| \leq 2^{m-l-1}$, $|c_{i,j}^m| \leq 2^{1-\delta_{ij}} \cdot 2^{2m-i-j-2} < 4^m$.

To show that Bott--Samelson rings are Koszul, we need the following definition of the complete intersection ring.

\begin{definition}
Let $A$ be a finitely generated algebra over a field $k$. Say $A$ is a \emph{complete intersection ring}, if $A\cong k[x_1,\cdots,x_n]/(f_1,\cdots, f_l)$, where the polynomial sequence $(f_1,\cdots,f_l)$ is a regular sequence in $k[x_1,\cdots,x_n]$. This means that $f_{1}\neq 0$ and $f_{i+1}$ is not a zero divisor in  $k[x_1,\cdots,x_n]/(f_1,\cdots,f_{i})$ for $1\leq i\leq l-1$.
\end{definition}

A finitely generated commutative algebra $A$ over a field $k$ is said to be a \emph{quadratic algebra} if $A\cong k[x_1,\cdots,x_n]/(f_1,\cdots,f_l)$ and all $f_i$ with $1\leq i \leq l$ are quadratic forms over $k$. Clearly, the Bott--Samelson ring $\overline{\mathrm{BS}}(\underline{w})$ is a quadratic algebra over $\mathbb{R}$.

Now we introduce the definition of Koszul algebra.

\begin{definition}
A \emph{Koszul algebra} $A$ over $k$ is a graded $k$-algebra such that the ground field $k$ has a linear (graded) free resolution, that is, there is an exact sequence

\begin{equation} 
    \cdots \rightarrow (A(-i))^{b_i} \rightarrow \cdots \rightarrow (A(-1))^{b_1} \rightarrow A \rightarrow k \rightarrow 0
\end{equation}
for some non-negative integers $b_{i}$. Recall that $A(-j)$ is the graded algebra $A$ with grading shifted up by $j$. The chain maps are all given by matrices with entries being 0 or linear forms after a choice of bases.
\end{definition}

A typical example of a Koszul algebra is a polynomial ring over a field, since the Koszul complex gives the minimal graded free resolution of the ground field.

It is well-known that any Koszul algebra is a quadratic algebra, see \cite[Proposition 1.2.3.]{beilinson1996koszul}. Conversely, a theorem of Tate \cite{tate1957homology} says that any quadratic complete intersection $k$-algebra is Koszul, see also \cite[Remark 1.12.] {conca2014koszul} for an easier argument.

Now we prove that any Bott--Samelson ring $\overline{\mathrm{BS}}(\underline{w})$ is Koszul. 

\begin{theorem} \label{thm-Koszul}
For any expression $\underline{w}=$ $\left(s_1, \ldots, s_n\right)$ in an arbitrary Coxeter system $(W, S)$, the Bott--Samelson ring $\overline{\mathrm{BS}}(\underline{w})$, defined in \eqref{def-gen}, is a Koszul $\mathbb{R}$-algebra.
\end{theorem}

\begin{proof}
    By Tate's Theorem, we only need to verify that, for an arbitrary $\underline{w}=$ $\left(s_1, \ldots, s_n\right)$ as above, $\overline{\mathrm{BS}}(\underline{w})$ is a complete intersection (since it is clear quadratic), that is, we need to show that $(q_1,\cdots,q_n)$ is a regular sequence in $\mathbb{R}[x_1,\cdots,x_n]$. We proceed from definition and only need to show that $q_{i+1}$ is not a zero divisor in $Q_i:=\mathbb{R}[x_1,\cdots,x_n]/(q_1,\cdots,q_{i})$.

    Note that $Q_i=(\mathbb{R}[x_1,\cdots,x_i]/(q_1,\cdots,q_i))[x_{i+1},\cdots,x_n]$, since $q_1,\cdots,q_i$ are polynomials with indeterminates $x_1,\cdots,x_{i}$. 
    In terms of expression \eqref{eq-expan}, we can write $q_{i+1}$ as $x_{i+1}^{2}+s_{i}$, with $s_{i}\in \mathbb{R}[x_1,\cdots,x_i]$. We order polynomials in $Q_i$ by graded (with respect to the total degree in $x_{i+1}, \cdots, x_n$) lexicographical order induced by $x_n > \cdots > x_{i+2} > x_{i+1}$.

    Denote the ideal $S_{i} := (q_1,\cdots,q_{i})\mathbb{R}[x_1,\cdots,x_n]$, which is the equivalent class of $0$ in $Q_i$. Assume, for the sake of contradiction, that $q_{i+1}$ is a zero divisor in $Q_i$. In other words,  there exists some $g \in Q_i$ but $g \notin S_{i}$ such that $gq_{i+1}\in S_{i}$. Write 
    \[g=\sum_{J=(j_1,\cdots,j_{n-i})}g_{J}x_{i+1}^{j_1}\cdots x_{n}^{j_{n-i}},\]
    where $g_{J}\in \mathbb{R}[x_1,\cdots,x_i]$ and we can assume all of them are not in $S_i$ by removing those in $S_i$ since $S_i$ is an ideal. Similarly, $gq_{i+1}$ is in the form of 
    $$\sum_{J=(j_1,\cdots,j_{n-i})}f_{J}x_{i+1}^{j_1}\cdots x_{n}^{j_{n-i}}$$
    with all $f_{J}\in \mathbb{R}[x_1,\cdots,x_i]$. One can easily deduce from $gq_{i+1} \in S_{i}$ that $f_J\in S_{i}$ for all $J$.

    Let $g_0 = g_{J_0}x_{i+1}^{j_{0,1}}\cdots x_{n}^{j_{0,n-i}}$ be the leading monomial of $g$ with respect to the above graded lexicographical order. According to the expression \eqref{eq-expan}, the leading monomial of $gq_{i+1}$ equals $g_{J_0}x_{i+1}^{j_{0,1}+2}\cdots x_{n}^{j_{0,n-i}}$. However, $g_{J_{0}}$ here, inherited from $g_0$, does not belong to $S_{i}$, which contradicts to the fact that $f_{J} \in S_i$ for all $J$. This contradiction implies that $q_{i+1}$ is not a zero divisor in $Q_i$. The proof is completed.
\end{proof}

\begin{eg}
    Consider an expression $\underline{w}=$ $\left(s_1, s_2, \cdots, s_n\right)$ in a Coxeter system $(W, S)$ such that all $m_{s_is_j}$'s are equal to $2$, that is, any $s_i$ and $s_j$ commute if $i \ne j$. Then the Bott--Samelson ring $\overline{\mathrm{BS}}(\underline{w})$ has the following presentation:
\begin{equation*} 
\overline{\mathrm{BS}}\left(s_1, s_2, \cdots, s_n\right) \cong \frac{\mathbb{R}\left[x_1, x_2, \ldots, x_n\right]}{\left(x_1^2, x_2^2, \cdots, x_n^2\right)}.
\end{equation*}
It is well known that this ring is Koszul, see, for example \cite[Example 6.3]{zbMATH07597716}.
\end{eg}

Inductively, using $q_1,\cdots,q_i$ to eliminate square terms $x_{1}^{2},\cdots,x_{i}^{2}$ in $q_{i+1}$ for all $i$, we can get a class of new generators of the ideal $(q_1,\cdots,q_n)$, denoted by $\{q_{1}',\cdots,q_{n}'\}$. We can express all of them as the following form:

\begin{equation} \label{redu-grb-basis}
    q_{k}' = x_k^2-\sum_{1 \leq i < j \leq k-1}d_{i,j}^{k}x_ix_j, \quad k=1, 2, \cdots, n
\end{equation}
where $d_{i,j}^{k}$'s are some real numbers whose explicit values are not needed in this paper.

We claim that the above $\{q_{1}',\cdots,q_{n}'\}$ gives a reduced Gr{\"o}bner basis of the ideal $(q_1,\cdots,q_n)$.

\begin{definition}
    Let $A=k[x_1,\cdots,x_n]$ be a polynomial ring over a field $k$ and fix a monomial ordering. Let $I$ be an ideal of $A$ and let $G$ be a finite generating set of $I$. The set $G$ is called a \emph{Gr{\"o}bner basis} of $I$ (with respect to the monomial ordering) if the ideal generated by the leading monomials of the polynomials in $I$ equals the ideal generated by the leading monomials of the polynomials in $G$. A Gröbner basis $G$ is \emph{reduced} if every leading monomial of the polynomials in $G$ is monic and does not divide any other monomials of the polynomials in $G$. 
\end{definition}

It is known that reduced Gröbner basis of an ideal $I$ (for a fixed monomial ordering) is unique \cite[Theorem 1.8.7.]{adams2022introduction} and the monomials in $A$ that are not divided by any leading monomials of the polynomials in the reduced Gröbner basis form a basis of the $k$-vector space $A / I$ \cite[[Proposition 2.1.6.]{adams2022introduction}. 

We have the following proposition, which gives another proof of the koszulness of the Bott--Samelson rings, since it is well known that a quadratic algebra is Koszul if the defining ideal of relations have a quadratic Gr{\"o}bner basis, see for example \cite[Section 4]{zbMATH01532175}. Note that a quadratic complete intersection may not have a quadratic Gr{\"o}bner basis with respect to any system of coordinates and any term order, see \cite{zbMATH00711325}.
\begin{proposition} \label{prop-Grobner}
    $\{q_{1}, \cdots, q_{n}\}$ is a Gr{\"o}bner basis and $\{q_{1}', \cdots, q_{n}'\}$ is the reduced Gr{\"o}bner basis of the ideal $(q_1,\cdots,q_n)$ with respect to the lexicographical order induced by $x_n > \cdots > x_2 > x_1$. As a corollary, $\displaystyle \frac{\mathbb{R}\left[x_1, x_2, \ldots, x_n\right]}{\left(q_1, q_2, \ldots, q_n\right)}$ has a basis $\{x_1^{\varepsilon_1}x_2^{\varepsilon_2}\cdots x_n^{\varepsilon_n} \mid \varepsilon_i \in \{0,1\} ,i \in [n]\}$.
\end{proposition}

\begin{proof}
 The leading monomial of $q_{k}$ as well as $q_{k}'$ is clearly $x_{k}^{2}$. $x_{k}^{2}$ and $x_{j}^{2}$ with $k\neq j$ have no common factors except constants. Hence by \cite[Lemma 5.66 and Theorem 5.68]{becker1993grobner}, the generating set $\{q_{1}, \cdots, q_{n}\}$ and $\{q_{1}',\cdots,q_{n}'\}$ are both Gr{\"o}bner bases.

    For the monic quadratic $x_{k}^2$ in $q_{k}'$, all other monomials of the polynomials in $\{q_{1}', \cdots, q_{n}'\}$ are either $x_{j}^{2}$ with $j\neq k$ or $d_{i,l}^{m}x_{i}x_{l}$ with $i\neq l$. So $\{q_{1}',\cdots,q_{n}'\}$ satisfies the definition of reduced Gr{\"o}bner basis.  
\end{proof}

\section{Hodge theory of Bott--Samelson Rings} \label{sec-Hodge}

In this section we prove that Bott--Samelson Rings satisfy properties in the K\"ahler package. We refer to \cite[Section 2]{elias2014hodge}, \cite[Chapter 17]{elias2020introduction}, and \cite[Section 2]{gui2024equivariant} for these properties in the linear algebra context.

For the expression $\underline{w}=(s_1,\cdots,s_n)$ in an arbitrary Coxeter system $(W,S)$, the Bott--Samelson ring $\overline{\mathrm{BS}}(\underline{w})$ is a graded ring of the form \eqref{def-gen}. It has a graded basis $\{x_1^{\varepsilon_1}x_2^{\varepsilon_2}\cdots x_n^{\varepsilon_n} \mid \varepsilon_i \in \{0,1\} ,i \in [n]\}$ when seen as a graded vector space over $\mathbb{R}$ (see Proposition \ref{prop-Grobner}). So given any degree $d$ homogeneous polynomial in $\overline{\mathrm{BS}}(\underline{w})$, it can be expressed uniquely as a real linear combination of some degree $d$ elements in the above basis of $\overline{\mathrm{BS}}(\underline{w})$. In particular, any monomial with degree greater than $n$ must be equal to $0$ in $\overline{\mathrm{BS}}(\underline{w})$. Note that the map $\overline{\mathrm{BS}}(s_{1},s_{2},\cdots, s_{k}) \hookrightarrow \overline{\mathrm{BS}}(\underline{w})$ induced by the natural map $\mathbb{R}[x_1,\cdots,x_k] \longrightarrow \overline{\mathrm{BS}}(\underline{w})$ is an embedding for any $k\leq n$. We will identify  $\overline{\mathrm{BS}}(s_{1},s_{2},\cdots,s_{k})$ with its image in $\overline{\mathrm{BS}}(\underline{w})$ under this embedding. 

We define the degree map $\deg \colon \overline{\mathrm{BS}}(\underline{w}) \rightarrow \mathbb{R}$ by sending any $\eta \in \overline{\mathrm{BS}}(\underline{w})$ to the coefficient of $x_{1}x_{2}\cdots x_{n}$ when expressing $\eta$ as an $\mathbb{R}$-linear combination of the above basis. We are ready to prove the Poincar\'{e} duality for $\overline{\mathrm{BS}}(\underline{w})$.

\begin{proof}[Proof of Theorem \ref{thm-Hodge}(1)]\label{prf-Hodge-PD}

For any $k\leq n / 2$, we denote $\langle \cdot,\cdot \rangle_n$ the bilinear pairing 
$$\overline{\mathrm{BS}}^k(\underline{w}) \times \overline{\mathrm{BS}}^{n-k}(\underline{w}) \longrightarrow \mathbb{R}$$ by associating a pair $(\eta_{1},\eta_{2})\in \overline{\mathrm{BS}}(\underline{w})$ to $\deg(\eta_{1}\eta_{2})$.

We want to show that $\langle \cdot ,\cdot \rangle_n$ is non-degenerate for any $k\leq n / 2$.

For an arbitrary subset $I$ of $[n]$, we define $x^{I}:=x_1^{\varepsilon_1}x_2^{\varepsilon_2}\cdots x_n^{\varepsilon_n}$, where $\varepsilon_{i}=1$ if $\varepsilon_{i}\in I$ and $\varepsilon_{i}=0$ if $\varepsilon_{i}\notin I$. Note that we can rewrite the basis defined above as $B_n := \{x^{I} \mid I\subset [n] \}$ and one can easily see that a basis of $\overline{\mathrm{BS}}^k(\underline{w})$ is $B_n^k := \{x^{I} \mid I\subset [n] \ |I|=k \}$ for $1\leq k \leq n$. 

For a pair of basis elements $(x^I,x^J)\in \overline{\mathrm{BS}}^k(\underline{w}) \times \overline{\mathrm{BS}}^{n-k}(\underline{w})$, we claim that deg($x^{I}x^{J}$)=1 if $J=[n] \backslash I$, deg($x^{I}x^{J}$)=0 if otherwise. This claim implies the pairing $\langle \cdot  , \cdot \rangle_n$ is non-degenerate for all $k\leq n / 2$, since we have dual basis under each pairing.

We prove the claim now. If $J=[n] \backslash I$, $x^{I}x^{J}$ is exactly $x_{1}x_{2}\cdots x_{n}$, so deg($x^{I}x^{J}$)=1 by definition. For the remaining case, 
suppose $t$ is the largest index such that $t\in I\cap J$, then $x^{I}x^{J}$ is of the form $x_{1}^{\varepsilon_{1}}x_{2}^{\varepsilon_{2}}\cdots x_{t}^{2} x_{t+1}^{\varepsilon_{t+1}}\cdots x_{n}^{\varepsilon_{n}}$, where $\varepsilon_{t+1},\cdots,\varepsilon_{n} \leq 1$ and $x_{1}^{\varepsilon_{1}}x_{2}^{\varepsilon_{2}}\cdots x_{t}^{2}$ is of degree $\geq t$ . We expand $x_{t}^{2}$ as a linear combination of $x_{i}x_{j}$'s with $i,j < t$ by \eqref{redu-grb-basis}. So $x_{1}^{\varepsilon_{1}}x_{2}^{\varepsilon_{2}}\cdots x_{t}^{2}$ can be written as a linear combination of monomials with respect to $x_{1}, \cdots, x_{t-1}$, and each of the monomial has degree $\geq t$. Hence $x_{1}^{\varepsilon_{1}}x_{2}^{\varepsilon_{2}}\cdots x_{t}^{2}$ is equal to $0$ in $\overline{\mathrm{BS}}(s_{1},s_{2},\cdots,s_{t-1})$. Identify  $\overline{\mathrm{BS}}(s_{1},s_{2},\cdots,s_{t-1})$ with its image under the embedding above, then $x_{1}^{\varepsilon_{1}}x_{2}^{\varepsilon_{2}}\cdots x_{t}^{2}=0$ also holds in $\overline{\mathrm{BS}}(\underline{w})$, which implies $x^{I}x^{J}=x_{1}^{\varepsilon_{1}}x_{2}^{\varepsilon_{2}}\cdots x_{t}^{2} x_{t+1}^{\varepsilon_{t+1}}\cdots x_{n}^{\varepsilon_{n}}=0\in \overline{\mathrm{BS}}(\underline{w})$. The claim is proven and Proof of Theorem \ref{thm-Hodge}(1) is completed.

\end{proof}

To complete the proof of the hard Lefschetz theorem (hL) and the Hodge--Riemann bilinear relations (HR), we need the following lemma (in our degree convention), which enhances \cite[Lemma 5.2]{elias2014hodge} and \cite[Lemma 18.31]{elias2020introduction}.
\begin{lemma}\label{lem-hrequiv}
    Let $V=\bigoplus_{i=0}^nV^i$ and $W=\bigoplus_{i=0}^{n-1}W^i$ be two finite-dimensional graded vector spaces, equipped
with graded non-degenerate forms and Lefschetz operators $(\langle \cdot,\cdot \rangle_V, L_V)$ and $(\langle \cdot,\cdot \rangle_W,L_W)$ (without assuming satisfying the hard Lefschetz theorem a priori as in \cite[Lemma 18.31]{elias2020introduction}). Assume that for each $i$, $\dim V^i = \dim W^{i}+\dim W^{i-1}$ and if for all $k \le n/2$, the signature of $(\cdot,\cdot)_{L_W}$ on $P^k(W):=\operatorname{ker}(L_{W}^{n-2k}\vert_{W^{k}})$ 
equals the signature of $(\cdot,\cdot)_{L_V}$ on all of $V^{k}$, 
then $\operatorname{HR}(V,\langle \cdot,\cdot \rangle_V,L_V)$ and $\operatorname{HR}(W,\langle \cdot,\cdot \rangle_W,L_W)$ are equivalent.
\end{lemma}

\begin{proof}
    Firstly, we show that $\operatorname{HR}(W,\langle \cdot,\cdot \rangle_W,L_W)$ implies $\operatorname{HR}(V,\langle \cdot,\cdot \rangle_V,L_V)$ under the above assumptions.

    We do induction on $i$ of $V^{i}$. When $i=0$, we have $\dim V^{0}=\dim W^{0}$, $P^{0}(W)=W^{0}$ and $P^{0}(V)=V^{0}$. By the assumption on the signatures and $\operatorname{HR}(W,\langle \cdot,\cdot \rangle_W,L_W)$, we know that $\dim V^{0}=\dim W^{0}=\operatorname{sign}(\cdot,\cdot)_{L_W}^{0}=\operatorname{sign}(\cdot,\cdot)_{L_V}^{0}$. This shows that $(\ ,\ )_{V}$ is positive definite on $V^{0}$.
    
    Suppose under the above assumptions, $\operatorname{HR}(W,\langle \cdot,\cdot \rangle_W,L_W)$ implies the Hodge--Riemann relations for $(V^i,\langle \cdot,\cdot \rangle_V,L_V)$ when $i\leq t-1$, thus implies $\operatorname{hL}(V^i,\langle \cdot,\cdot \rangle_V,L_V)$ for $i\leq t-1$. In particular, $L_{V}^{n-2t+2}: V^{t-1}\rightarrow V^{n-t+1}$ is bijective and $L_{V}: V^{t-1}\rightarrow V^{t}$ is injective. This leads to a decomposition $V^{t}=L_{V}(V^{t-1})\oplus P^{t}(V)$, where $P^{t}(V):=\operatorname{ker}(L_{V}^{n-2t+1}\vert_{V^{t}})$ and $L_{V}(V^{t-1})$ are orthogonal under $(\cdot,\cdot)_{L_V}^{t}$. Now we have,
    \begin{equation*}
    \begin{aligned}
&\operatorname{sign}(\cdot,\cdot)_{L_V}^{t}\vert_{P^{t}(V)}\\=&\operatorname{sign}(\cdot,\cdot)_{L_V}^{t}-\operatorname{sign}(\cdot,\cdot)_{L_V}^{t-1}\\=&\operatorname{sign}(\cdot,\cdot)_{L_W}^{t}\vert_{P^{t}(W)}-\operatorname{sign}(\cdot,\cdot)_{L_W}^{t-1}\vert_{P^{t-1}(W)}\\=&(-1)^{t}(\operatorname{dim}W^{t}-\operatorname{dim}W^{t-1})-(-1)^{t-1}(\operatorname{dim}W^{t-1}-\operatorname{dim}W^{t-2})\\=&(-1)^{t}(\operatorname{dim}V^{t}-\operatorname{dim}V^{t-1})\\=&(-1)^{t}(\operatorname{dim}P^{t}(V)),
\end{aligned}\end{equation*}
where the first and the last equalities follow from the orthogonal direct sum decomposition, while the other equalities follow from our assumptions. Therefore, $(\cdot,\cdot)_{L_V}\vert_{P^{t}(V)}$ is $(-1)^{t}$-definite, that is, Hodge--Riemann bilinear relations hold for $V^{t}$.

    The arguments for the reverse implication are similar by induction on $i$ of $W^{i}$ with the above equations replaced by

     \begin{equation*}
    \begin{aligned}
&\operatorname{sign}(\cdot,\cdot)_{L_W}^{t}\vert_{P^{t}(W)}\\=&\operatorname{sign}(\cdot,\cdot)_{L_V}^{t}\\=&\sum_{k=0}^{t}(-1)^{k}(\dim V^{k}-\dim V^{k-1})
    \\=&\sum_{k=0}^{t}(-1)^{k}(\dim W^{k}-\dim W^{k-2})\\=&(-1)^{t}(\dim W^{t}-\dim W^{t-1})\\=&(-1)^{t}\dim P^{t}(W)
    \end{aligned}
    \end{equation*}
We leave the details to the reader. The proof is completed.
\end{proof}

\begin{proof}[Proof of \ref{thm-Hodge}(2)(3)]\label{prf-hlhr}
    We use induction on $n = \ell(\underline{w})$. One directly checks that ($\mathrm{hL}$) and ($\mathrm{HR}$) are satisfied for any expression $\underline{w} = (s_1)$ in an arbitrary Coxeter system $(W,S)$ whenever $c_1 > 0$. Fix any expression $\underline{w} = (s_1,\ldots,s_n)$ in a Coxeter system $(W,S)$ and $k \leq n / 2$, we denote $(s_1,\ldots,s_{n-1})$ by $\underline{w}'$. By induction, there exists $\{c_i\}_{i = 1}^{n-1}$ with $c_i>0$, such that ($\mathrm{hL}$) and ($\mathrm{HR}$) hold in $\overline{\mathrm{BS}}(\underline{w}')$ with respect to $\ell_{n-1}=\sum_{1\leq i\leq n-1} c_i x_i$.

Since we have proved the Poincaré duality for $\overline{\mathrm{BS}}(\underline{w})$, the Hodge--Riemann relations for $\overline{\mathrm{BS}}(\underline{w})$ are stronger and imply the hard Lefschetz theorem for $\overline{\mathrm{BS}}(\underline{w})$. By Lemma \ref{lem-hrequiv} and induction, it suffices to show that there exists some $c_n>0$, such that for all $k\le n / 2$, $(\cdot,\cdot)_{\ell_{n-1}}$ on ${P_{\ell_{n-1}}^{k}}$ has the same signature as $(\cdot,\cdot)_{\ell_n}$ on $\overline{\mathrm{BS}}^k(\underline{w})$.
    
    If $k = n/2$, then $(\cdot,\cdot)_{\ell_{n-1}}$ has signature $0$ on ${P_{\ell_{n-1}}^{k}}$ since ${P_{\ell_{n-1}}^{k}}=0$ for $\overline{\mathrm{BS}}(\underline{w}')$, and the Gram matrix of $(\cdot,\cdot)_{\ell_n}$ (which is actually the Poincaré pairing $\langle \cdot ,\cdot \rangle_n$ and is independent of $\ell_n$) on $\overline{\mathrm{BS}}^k(\underline{w})$ with respect to the basis $B_n^k$ is \[\begin{pmatrix}
         & & & 1\\ & & \begin{turn}{90}$\ddots$\end{turn} & \\ & 1 & & \\ 1 & & & 
    \end{pmatrix}\] 
    (by orthogonality of the basis), which is of signature $0$.

    Now suppose $k < n/2$, then we have 
    \[\ell_n^{n-2k} = \left(\ell_{n-1}+c_n x_n\right)^{n-2k} = \ell_{n-1}^{n-2k}+(n-2k)c_nx_n\ell_{n-1}^{n-2k-1}+c_n^2\Tilde{\ell}_n,\] where $\Tilde{\ell}_n \in \overline{\mathrm{BS}}^{n-2k}(\underline{w})$ depends polynomially on $c_n$. 
    We decompose $\overline{\mathrm{BS}}^k(\underline{w})$ into the direct sum of three subspaces: $P_{\ell_{n-1}}^{k}$, $\ell_{n-1}(\overline{\mathrm{BS}}^{k-1}(\underline{w}'))$ and $\overline{\mathrm{BS}}^{k-1}(\underline{w}')x_n$. That is, 
    \[\overline{\mathrm{BS}}^k(\underline{w}) = P_{\ell_{n-1}}^{k} \oplus \ell_{n-1}(\overline{\mathrm{BS}}^{k-1}(\underline{w}')) \oplus \overline{\mathrm{BS}}^{k-1}(\underline{w}')x_n.\]
    Choose bases $B_1'$ and $B_2'$ of $P_{\ell_{n-1}}^{k}$ and $\overline{\mathrm{BS}}^{k-1}(\underline{w}')$, respectively. Then $B' = B_1' \cup \ell_{n-1}B_2' \cup B_2'x_n$ forms a basis of $\overline{\mathrm{BS}}^k(\underline{w})$. Here $\ell_{n-1}B_2' := \{\ell_{n-1}\eta \mid \eta \in B_2'\}$ and $B_2'x_n$ is similarly defined. 

    Now consider the bilinear form $(\cdot,\cdot)_{\ell_n}$ on $\overline{\mathrm{BS}}^k(\underline{w})$. 
    
    For $\eta_1,\eta_2 \in B_1'$, we have 
    \begin{equation}\label{eq-bf-pp}
        \begin{split}
            (\eta_1,\eta_2)_{\ell_n}&= \left\langle \left(\ell_{n-1}^{n-2k}+(n-2k)c_nx_n\ell_{n-1}^{n-2k-1}+c_n^2\Tilde{\ell}_n\right)\eta_1, \eta_2\right\rangle_n 
            \\ & = (n-2k)c_n\langle \ell_{n-1}^{n-2k-1}\eta_1 \cdot x_n,\eta_2\rangle_n+ c_n^2 \langle \Tilde{\ell}_n\eta_1,\eta_2\rangle_n
            \\ & = (n-2k)c_n\langle \ell_{n-1}^{n-2k-1}\eta_1,\eta_2\rangle_{n-1}+ c_n^2 \langle \Tilde{\ell}_n\eta_1,\eta_2\rangle_n 
            \\ & = (n-2k)c_n\left( \eta_1,\eta_2\right)_{\ell_{n-1}}+ c_n^2 \langle \Tilde{\ell}_n\eta_1,\eta_2\rangle_n,
        \end{split}
    \end{equation}
the second equation follows from $\ell_{n-1}^{n-2k} \eta_1= 0$ since $\eta_1 \in P_{\ell_{n-1}}^{k}$. 

For $\eta_1 \in B_1'$ and $\eta_2 \in B_2'$, we have 
    \begin{equation}\label{eq-bf-pl}
        \begin{split}
            (\eta_1,\ell_{n-1}\eta_2)_{\ell_n} &= \left\langle \left(\ell_{n-1}^{n-2k}+(n-2k)c_nx_n\ell_{n-1}^{n-2k-1}+c_n^2\Tilde{\ell}_n\right)\eta_1, \ell_{n-1}\eta_2\right\rangle_n 
            \\ & = (n-2k)c_n\langle \ell_{n-1}^{n-2k}\eta_1 \cdot x_n,\eta_2\rangle_n+ c_n^2 \langle \Tilde{\ell}_n\eta_1,\ell_{n-1}\eta_2\rangle_n
            \\ & = c_n^2 \langle \Tilde{\ell}_n\eta_1,\ell_{n-1}\eta_2\rangle_n
        \end{split}
    \end{equation}
    and
    \begin{equation}\label{eq-bf-pn}
        \begin{split}
        (\eta_1,\eta_2x_n)_{\ell_n}&= \left\langle \left(\ell_{n-1}^{n-2k}+(n-2k)c_nx_n\ell_{n-1}^{n-2k-1}+c_n^2\Tilde{\ell}_n\right)\eta_1, \eta_2x_n\right\rangle_n
            \\ & = (n-2k)c_n\langle \ell_{n-1}^{n-2k-1}\eta_1\cdot x_n,\eta_2\cdot x_n\rangle_n+ c_n^2 \langle \Tilde{\ell}_n\eta_1,\eta_2x_n\rangle_n
            \\ & = c_n^2 \langle \Tilde{\ell}_n\eta_1,\eta_2x_n\rangle_n.
        \end{split}
    \end{equation}
    We note that the third equation of \eqref{eq-bf-pn} holds due to the orthogonal basis of the Poincar\'e pairing and that both $\ell_{n-1}^{n-2k-1}\eta_1x_n$ and $\eta_2x_n$ contain $x_n$.

    For $\eta_1,\eta_2 \in B_2'$, we have 
    \begin{equation}\label{eq-bf-ll}
        \begin{split}
            (\ell_{n-1}\eta_1,\ell_{n-1}\eta_2)_{\ell_n} &= \left\langle \left(\ell_{n-1}^{n-2k+1}+(n-2k)c_nx_n\ell_{n-1}^{n-2k}+c_n^2\Tilde{\ell}_n\ell_{n-1}\right)\eta_1, \ell_{n-1}\eta_2\right\rangle_n 
            \\ & = (n-2k)c_n\langle \ell_{n-1}^{n-2k}\eta_1 \cdot x_n,\ell_{n-1}\eta_2\rangle_n+ c_n^2 \langle \Tilde{\ell}_n\ell_{n-1}\eta_1,\ell_{n-1}\eta_2\rangle_n
            \\ & = (n-2k)c_n\langle \ell_{n-1}^{n-2k+1}\eta_1 ,\eta_2\rangle_{n-1}+ c_n^2 \langle \Tilde{\ell}_n\ell_{n-1}\eta_1,\ell_{n-1}\eta_2\rangle_n
            \\ & = (n-2k)c_n\left(\eta_1 ,\eta_2\right)_{\ell_{n-1}}+ c_n^2 \langle \Tilde{\ell}_n\ell_{n-1}\eta_1,\ell_{n-1}\eta_2\rangle_n.
        \end{split}
    \end{equation}
    Note that $\left\langle \ell_{n-1}^{n-2k+1}\eta_1, \ell_{n-1}\eta_2\right\rangle_n=0$, since neither $\ell_{n-1}^{n-2k+1}\eta_1$ nor $\ell_{n-1}\eta_2$ contains $x_n$. Similarly, 
    \begin{equation}\label{eq-bf-ln}
        \begin{split}
            (\ell_{n-1}\eta_1,\eta_2x_n)_{\ell_n}&= \left\langle \left(\ell_{n-1}^{n-2k+1}+(n-2k)c_nx_n\ell_{n-1}^{n-2k}+c_n^2\Tilde{\ell}_n\ell_{n-1}\right)\eta_1, \eta_2x_n\right\rangle_n 
            \\ & = \langle\ell_{n-1}^{n-2k+1}\eta_1,\eta_2x_n\rangle_n + c_n^2 \langle \Tilde{\ell}_n\ell_{n-1}\eta_1,\eta_2x_n\rangle_n
            \\ & = \langle\ell_{n-1}^{n-2k+1}\eta_1,\eta_2\rangle_{n-1} + c_n^2 \langle \Tilde{\ell}_n\ell_{n-1}\eta_1,\eta_2x_n\rangle_n
            \\ & = \left(\eta_1 ,\eta_2\right)_{\ell_{n-1}}+ c_n^2 \langle \Tilde{\ell}_n\ell_{n-1}\eta_1,\eta_2x_n\rangle_n.
        \end{split}
    \end{equation}
    Finally,
    \begin{equation}\label{eq-bf-nn}
        \begin{split}
    (\eta_1x_n,\eta_2x_n)_{\ell_n}&= \left\langle \left(\ell_{n-1}^{n-2k}+(n-2k)c_nx_n\ell_{n-1}^{n-2k-1}+c_n^2\Tilde{\ell}_n\right)\eta_1x_n, \eta_2x_n\right\rangle_n 
            \\ & = (n-2k)c_n\langle\ell_{n-1}^{n-2k-1}x_n^2\eta_1,\eta_2x_n\rangle_n + c_n^2 \langle \Tilde{\ell}_n\eta_1x_n,\eta_2x_n\rangle_n
            \\ & = (n-2k)c_n\left\langle\sum_{1 \leq i \leq j \leq n-1} c_{i,j}^n x_ix_j\ell_{n-1}^{n-2k-1}\eta_1,\eta_2\right\rangle_{n-1} + c_n^2 \langle \Tilde{\ell}_n\eta_1x_n,\eta_2x_n\rangle_n,
        \end{split}
    \end{equation}
    where $c_{i,j}^n$ is defined in \eqref{eq-expan}.
    
    From \eqref{eq-bf-pp} to \eqref{eq-bf-nn} above, we deduce that the Gram matrix of $(\cdot,\cdot)_{\ell_n}$ on $\overline{\mathrm{BS}}^k(\underline{w})$ with respect to the basis $B'$ has the following form:
    \begin{equation}\label{eq-gram}
        \begin{pmatrix}
    (n-2k)c_nG_{11} + c_n^2\Tilde{G}_{11}& c_n^2\Tilde{G}_{12}& c_n^2\Tilde{G}_{13}\\
    c_n^2\Tilde{G}_{21}& (n-2k)c_nG_{22}+c_n^2\Tilde{G}_{22} & G_{23}+c_n^2\Tilde{G}_{23}\\
    c_n^2\Tilde{G}_{31} & G_{32} +c_n^2\Tilde{G}_{32}& (n-2k)c_nG_{33}+c_n^2\Tilde{G}_{33} \\
\end{pmatrix}.
\end{equation}
    
The meanings of $G_{ij}$'s and $\Tilde{G}_{ij}$'s are as follows: $G_{11}$ is the Gram matrix of $\left(\cdot,\cdot\right)_{\ell_{n-1}}$ on $P_{\ell_{n-1}}^{k}$ with respect to the basis $B_1'$. $G_{22} = G_{23} = G_{32}$ is the Gram matrix of $(\cdot,\cdot)_{\ell_{n-1}}$ with respect to the basis $B_2'$. $G_{33}$ is the Gram matrix of \[\left\langle\left(\sum_{1 \leq i \leq j \leq n-1} \ell_{n-1}^{n-2k-1}c_{i,j}^n x_ix_j\right)\cdot,\cdot\right\rangle_{n-1} \colon \overline{\mathrm{BS}}^{k-1}(\underline{w}') \times \overline{\mathrm{BS}}^{k-1}(\underline{w}') \longrightarrow \mathbb{R}\] with respect to the basis $B_2'$. $(\Tilde{G}_{ij})$ is the Gram matrix of $\left(\cdot,\cdot\right)_{\Tilde{\ell}_n}$ on $\overline{\mathrm{BS}}^k(\underline{w}) $
with respect to the basis $B'$. 

By construction, $\Tilde{G}_{ij}$ is bounded as $c_n \to 0$. By induction, $G_{11}$ is $(-1)^k$-definite. Hence so is $(n-2k)c_nG_{11}+c_n^2\Tilde{G}_{11}$ as $c_n \to 0$. We reduce the matrix displayed in \eqref{eq-gram} to the following form after congruent transformations when $c_n > 0$ is sufficiently small:
$$\begin{pmatrix}
    (n-2k)c_nG_{11} + c_n^2\Tilde{G}_{11}& O& O\\
    O & (n-2k)c_nG_{22}+c_n^2M_{22} & G_{23}+c_n^2M_{23}\\
    O & G_{32} +c_n^2M_{32}& (n-2k)c_nG_{33}+c_n^2M_{33} \\
\end{pmatrix}.
$$

It is not hard to see that $M_{ij}$'s are also bounded as $c_n \to 0$ by looking at the order of $c_n$. Let $\operatorname{sgn}(\cdot)$ be the signature of a given matrix. When $c_n>0$ is sufficiently small, we have
\[\operatorname{sgn}\begin{pmatrix}
    (n-2k)c_nG_{22}+c_n^2M_{22} & G_{23}+c_n^2M_{23}\\
    G_{32} +c_n^2M_{32}& (n-2k)c_nG_{33}+c_n^2M_{33}\\
\end{pmatrix} = \operatorname{sgn} \begin{pmatrix}
    O & G_{23}\\
    G_{32}& O\\
\end{pmatrix} = 0,\] since by induction $G_{23} = G_{32}$ is an invertible symmetric matrix. Here $\operatorname{sgn} \begin{pmatrix}
    O & G_{23}\\
    G_{32}& O\\
\end{pmatrix} = 0$ is from basic linear algebra, since we can assume $G_{23} = G_{32}$ is diagonal after changing of basis.

Also as $c_n \to 0$, the following identities of signatures hold:
\begin{equation*}
    \begin{split}
        & \operatorname{sgn}\begin{pmatrix}
    (n-2k)c_nG_{11} + c_n^2\Tilde{G}_{11}& O& O\\
    O & (n-2k)c_nG_{22}+c_n^2M_{22} & G_{23}+c_n^2M_{23}\\
    O & G_{32} +c_n^2M_{32}& (n-2k)c_nG_{33}+c_n^2M_{33} \\
\end{pmatrix} \\ =& \operatorname{sgn} \left((n-2k)c_nG_{11} + c_n^2\Tilde{G}_{11}\right)\\ &+ \operatorname{sgn} \begin{pmatrix}
    (n-2k)c_nG_{22}+c_n^2M_{22} & G_{23}+c_n^2M_{23}\\
    G_{32} +c_n^2M_{32}& (n-2k)c_nG_{33}+c_n^2M_{33}\\
\end{pmatrix}  \\= &\operatorname{sgn} G_{11},
    \end{split}
\end{equation*}
which complete the proof of that $(\cdot,\cdot)_{\ell_{n-1}}$ on ${P_{\ell_{n-1}}^{k}}$ has the same signature as $(\cdot,\cdot)_{\ell_n}$ on $\overline{\mathrm{BS}}^k(\underline{w})$. The proof of \ref{thm-Hodge}(2)(3) is completed.
\end{proof}

\section{Questions and Future work} \label{sec-question}
We would like to finish with some questions and conjectures. Since the K\"ahler package and Koszulness are quite strong restrictions of a graded ring, we would like to ask the following
\begin{question}
Do the K\"ahler package and Koszulness of the Bott--Samelson rings have some implications for the study of Coxeter systems?
\end{question}

In \cite{richmond2021isomorphism}, the authors studied the isomorphism problem for Schubert varieties in the full flag variety of Kac--Moody type and gave a beautiful necessary and sufficient 
criterion for when two
such Schubert varieties (from potentially different flag varieties) are isomorphic, in terms
of the Cartan matrix and reduced words for the indexing Weyl group elements. Since we prove that the cohomology ring of the Bott--Samelson variety is a combinatorial invariant which only depends on the word in the Weyl group, it is natural to ask the following 
\begin{question} \label{ques-iso}
When are two Bott--Samelson varieties (associated with words from potentially different Weyl groups) isomorphic?
\end{question}

In \cite{zbMATH06579100}, the author studied a special fiber of the map of some Bott--Samelson variety into the flag variety in general types, known today as the ``brick manifolds''. She showed that a subfamily of brick manifolds gives a resolution of singularities for (closed) Richardson varieties $X_u^v$, which are certain (usually) singular varieties formed by the intersection of Schubert varieties and opposite Schubert varieties within a flag variety. For $u=e$ and $v$ a Coxeter element in type A, it is showed in \cite{zbMATH07094320} that the corresponding brick manifold is a non-commutative analogue of the compactified Deligne–Mumford moduli
spaces $\overline{\mathcal{M}}_{0, n+1}$ and also can be viewed as a projective wonderful
model of a certain hyperplane arrangement with respect to its minimal building set in the sense of de Concini and Procesi \cite{de1995wonderful}. Therefore, the cohomology ring of the brick manifold in this special case is isomorphic to the the Chow ring of the intersection lattice of the hyperplane arrangement with respect to the minimal building set in the sense of 
Feichtner and Yuzvinsky \cite{feichtner2004chow}, thus it satisfies the K\"ahler package. Also it has a quadratic Gr{\"o}bner basis, and therefore is Koszul by a result of Coron on supersolvable lattices \cite{coron2023supersolvabilitybuiltlatticeskoszulness}. It is reasonable to hope that some analogs of the results in this paper hold in this more general setting.
\begin{question} \label{ques-brick!}
Can one compute a presentation of the cohomology rings of brick manifolds? Are they combinatorial invariants? Do the results of K\"ahler package and Koszulness of the Bott--Samelson rings hold for these ``brick rings''?
\end{question}

\bibliographystyle{amsplain}
\bibliography{template}

\providecommand{\bysame}{\leavevmode\hbox to3em{\hrulefill}\thinspace}
\providecommand{\MR}{\relax\ifhmode\unskip\space\fi MR }
\providecommand{\MRhref}[2]{%
  \href{http://www.ams.org/mathscinet-getitem?mr=#1}{#2}
}
\providecommand{\href}[2]{#2}
\begin{thebibliography}{10}

\bibitem{adams2022introduction}
William~W. Adams and Philippe Loustaunau, \emph{An introduction to {Gr{\"o}bner} bases}, Grad. Stud. Math., vol.~3, Providence, RI: American Mathematical Society, 1994 (English).

\bibitem{adiprasito2018hodge}
Karim Adiprasito, June Huh, and Eric Katz, \emph{Hodge theory for combinatorial geometries}, Annals of Mathematics \textbf{188} (2018), no.~2, 381--452.

\bibitem{becker1993grobner}
Thomas Becker and Volker Weispfenning, \emph{Gr{\"o}bner bases: a computational approach to commutative algebra. {In} cooperation with {Heinz} {Kredel}}, Grad. Texts Math., vol. 141, New York: Springer-Verlag, 1993 (English).

\bibitem{beilinson2018faisceaux}
Alexander Beilinson, Joseph Bernstein, Pierre Deligne, and Ofer Gabber, \emph{Faisceaux pervers. {Actes} du colloque ``{Analyse} et {Topologie} sur les {Espaces} {Singuliers}''. {Partie} {I}}, 2nd edition ed., Ast{\'e}risque, vol. 100, Paris: Soci{\'e}t{\'e} Math{\'e}matique de France (SMF), 2018 (French).

\bibitem{beilinson1996koszul}
Alexander Beilinson, Victor Ginzburg, and Wolfgang Soergel, \emph{Koszul duality patterns in representation theory}, J. Am. Math. Soc. \textbf{9} (1996), no.~2, 473--527 (English).

\bibitem{bott1955cohomology}
Raoul Bott and Hans Samelson, \emph{The cohomology ring of {G/T}}, Proceedings of the National Academy of Sciences \textbf{41} (1955), no.~7, 490--493.

\bibitem{bott1958applications}
\bysame, \emph{Applications of the theory of {M}orse to symmetric spaces}, American Journal of Mathematics \textbf{80} (1958), no.~4, 964--1029.

\bibitem{conca2014koszul}
Aldo Conca, Sandra Di~Rocco, Jan Draisma, June Huh, Bernd Sturmfels, Filippo Viviani, and Aldo Conca, \emph{Koszul algebras and their syzygies}, Combinatorial Algebraic Geometry: Levico Terme, Italy 2013, Editors: Sandra Di Rocco, Bernd Sturmfels (2014), 1--31.

\bibitem{coron2023supersolvabilitybuiltlatticeskoszulness}
Basile Coron, \emph{Supersolvability of built lattices and {K}oszulness of generalized {C}how rings}, preprint, arXiv:2302.13072, 2023.

\bibitem{de1995wonderful}
C.~De~Concini and C.~Procesi, \emph{Wonderful models of subspace arrangements}, Sel. Math., New Ser. \textbf{1} (1995), no.~3, 459--494 (English).

\bibitem{demazure1974desingularisation}
Michel Demazure, \emph{Desingularisation des vari{\'e}t{\'e}s de {Schubert} g{\'e}n{\'e}ralis{\'e}es}, Ann. Sci. {\'E}c. Norm. Sup{\'e}r. (4) \textbf{7} (1974), 53--88 (French).

\bibitem{dotsenko2022homotopy}
Vladimir Dotsenko, \emph{Homotopy invariants for {{\(\overline{\mathcal{M}}_{0,n}\)}} via {Koszul} duality}, Invent. Math. \textbf{228} (2022), no.~1, 77--106 (English).

\bibitem{zbMATH07094320}
Vladimir Dotsenko, Sergey Shadrin, and Bruno Vallette, \emph{Toric varieties of {Loday}'s associahedra and noncommutative cohomological field theories}, J. Topol. \textbf{12} (2019), no.~2, 463--535 (English).

\bibitem{zbMATH00711325}
David Eisenbud, Alyson Reeves, and Burt Totaro, \emph{Initial ideals, {Veronese} subrings, and rates of algebras}, Adv. Math. \textbf{109} (1994), no.~2, 168--187 (English).

\bibitem{elias2020introduction}
Ben Elias, Shotaro Makisumi, Ulrich Thiel, and Geordie Williamson, \emph{Introduction to {Soergel} bimodules}, RSME Springer Ser., vol.~5, Cham: Springer, 2020 (English).

\bibitem{elias2014hodge}
Ben Elias and Geordie Williamson, \emph{The {Hodge} theory of {Soergel} bimodules.}, Ann. Math. (2) \textbf{180} (2014), no.~3, 1089--1136 (English).

\bibitem{elias2016kazhdan}
\bysame, \emph{Kazhdan-{Lusztig} conjectures and shadows of {Hodge} theory}, Arbeitstagung Bonn 2013. In memory of Friedrich Hirzebruch. Proceedings of the meeting, Bonn, Germany, May, 22--28, 2013, Basel: Birkh{\"a}user/Springer, 2016, pp.~105--126 (English).

\bibitem{elias2016soergel}
\bysame, \emph{Soergel calculus}, Represent. Theory \textbf{20} (2016), 295--374 (English).

\bibitem{zbMATH06579100}
Laura Escobar, \emph{Brick manifolds and toric varieties of brick polytopes}, Electron. J. Comb. \textbf{23} (2016), no.~2, research paper p2.25, 18 (English).

\bibitem{zbMATH07597716}
Eleonore Faber, Martina Juhnke-Kubitzke, Haydee Lindo, Claudia Miller, Rebecca Rebhuhn-Glanz, and Alexandra Seceleanu, \emph{Canonical resolutions over {Koszul} algebras}, Women in commutative algebra. Proceedings of the 2019 WICA workshop, Banff, Alberta, Canada, October 20--25, 2019, Cham: Springer, 2022, pp.~281--301 (English).

\bibitem{feichtner2004chow}
Eva~Maria Feichtner and Sergey Yuzvinsky, \emph{Chow rings of toric varieties defined by atomic lattices}, Inventiones mathematicae \textbf{155} (2004), no.~3, 515--536.

\bibitem{zbMATH01532175}
R.~Fr{\"o}berg, \emph{Koszul algebras}, Advances in commutative ring theory. Proceedings of the 3rd international conference, Fez, Morocco, New York, NY: Marcel Dekker, 1999, pp.~337--350 (English).

\bibitem{gui2024equivariant}
Tao Gui and Rui Xiong, \emph{Equivariant log-concavity and equivariant {K{\"a}hler} packages}, J. Algebra \textbf{657} (2024), 379--401 (English).

\bibitem{hansen1973cycles}
Hans~Christian Hansen, \emph{On cycles in flag manifolds}, Mathematica Scandinavica \textbf{33} (1973), no.~2, 269--274.

\bibitem{kazhdan1979representations}
David Kazhdan and George Lusztig, \emph{Representations of {C}oxeter groups and {H}ecke algebras}, Inventiones mathematicae \textbf{53} (1979), no.~2, 165--184.

\bibitem{mastroeni2023chow}
Matthew Mastroeni and Jason McCullough, \emph{Chow rings of matroids are {K}oszul}, Mathematische Annalen \textbf{387} (2023), no.~3, 1819--1851.

\bibitem{zbMATH07303202}
Chris McDaniel and Larry Smith, \emph{Equivariant coinvariant rings, {Bott}-{Samelson} rings and {Watanabe}'s bold conjecture}, J. Pure Appl. Algebra \textbf{225} (2021), no.~5, 50 (English), Id/No 106524.

\bibitem{richmond2021isomorphism}
Edward Richmond and William Slofstra, \emph{The isomorphism problem for {S}chubert varieties}, preprint, arXiv:2103.08114, 2021.

\bibitem{soergel1990kategorie}
Wolfgang Soergel, \emph{{Kategorie O} , perverse garben und moduln über den koinvarianten zur {W}eyl gruppe}, J. Am. Math. Soc. \textbf{3} (1990), no.~2, 421--445 (German).

\bibitem{soergel1992combinatorics}
\bysame, \emph{The combinatorics of {Harish}-{Chandra} bimodules}, J. Reine Angew. Math. \textbf{429} (1992), 49--74 (English).

\bibitem{zbMATH01131858}
\bysame, \emph{Kazhdan-{Lusztig} polynomials and a combinatorics for tilting modules}, Represent. Theory \textbf{1} (1997), 83--114 (English).

\bibitem{soergel2007kazhdan}
\bysame, \emph{{Kazhdan-Lusztig}-polynome und unzerlegbare bimoduln über polynomringen}, J. Inst. Math. Jussieu \textbf{6} (2007), no.~3, 501--525 (German).

\bibitem{tate1957homology}
John Tate, \emph{Homology of {N}oetherian rings and local rings}, Illinois Journal of Mathematics \textbf{1} (1957), no.~1, 14--27.

\end{thebibliography}

\end{document}